\documentclass[twocolumn]{jsiamletters}

\makeatletter
\def\Hline{%
\noalign{\ifnum0=`}\fi\hrule \@height 0.8pt \futurelet
\reserved@a\@xhline}
\makeatother

\usepackage{comment}

\usepackage{color}
\usepackage{ulem}

\usepackage{color}
\usepackage{tabularx}

\group{
Geometric Shape Generation
}
\affiliation{Institute of Mathematics for Industry, Kyushu University}{744, Motooka, Nishi-ku, Fukuoka, 819-0395, Japan.}
\affiliation{Research and Education Center for Comprehensive Science, Akita Prefectural University}{
84-4, Aza Ebinokuchi, Tsuchiya, Yurihonjo, 
Akita 015-0055, Japan}

\authorinfo{Miyuki Koiso}{1*}{koiso@imi.kyushu-u.ac.jp}
\authorinfo{Umpei Miyamoto}{2}{umpei@akita-pu.ac.jp}

\email{koiso@imi.kyushu-u.ac.jp}
\title{
Stability of hypersurfaces of constant mean curvature with free boundary in two parallel hyperplanes}
\abstract{
Surfaces with constant mean curvature (CMC) are critical points of the area with volume constraint. They serve as a mathematical model of surfaces of soap bubbles and tiny liquid drops. CMC surfaces are said to be stable if the second variation of the area is nonnegative for all volume-preserving variations satisfying the given boundary condition. In this paper, we examine the stability of CMC hypersurfaces in general Euclidean space possibly having boundaries on two parallel hyperplanes. We reveal the stability of equilibrium hypersurfaces without self-intersection for the first time in all dimensions. The analysis is assisted by numerical computations.
}
\begin{document}

\maketitle

\section{Introduction}\label{intro}
Surfaces with constant mean curvature (CMC) are critical points of a variational problem to extremize the area with volume constraint. They serve as a mathematical model of surfaces of soap bubbles and tiny liquid drops, for which we can ignore the gravitational energy. CMC surfaces are said to be stable if the second variation of the area is nonnegative for all volume-preserving variations satisfying the given boundary condition. Surfaces with minimal area amongst all nearby surfaces within the class of surfaces enclosing the same volume are stable CMC surfaces. Since only such surfaces are realized as natural phenomena, it is important to judge the stability for a CMC surface, which is difficult in general. Hence, 
in this paper we study the stability for one of the simplest boundary conditions. 
Since above concepts are generalized to higher dimensions, we study CMC hypersurfaces in the $(n+2)$-dimensional Euclidean space ${\mathbb R}^{n+2}$ ($n \ge 1$) having boundaries on two parallel hyperplanes.  

Let us call one of the coordinate axes in ${\mathbb R}^{n+2}$ the $z$-axis. 
Consider two parallel hyperplanes $\Pi_1:=\{z=z_1\}$, $\Pi_2:=\{z=z_2\}$ ($z_2>z_1$). Denote by $\Omega$ the closed domain bounded by $\Pi_1$, $\Pi_2$, and by $\Omega^\circ$ the interior of $\Omega$. 
Let ${\mathcal S}(\Omega)$be the set of all compact connected oriented smooth hypersurfaces in $\Omega$ possibly with free boundary in $\Pi_1\cup\Pi_2$. 
For any hypersurface $X \in {\mathcal S}(\Omega)$, we denote by $A(X)$ the $(n+1)$-dimensional volume of $X$, which we will call the ``area of $X$''. 
We denote by $V(X)$ the $(n+2)$-dimensional volume of the compact domain $G(X)$ bounded by $X\cup\Pi_1\cup\Pi_2$, which we will call the ``volume of $X$''. 
For $V>0$, denote by ${\mathcal S}(\Omega; V)$ the set of all hypersurfaces $X$ in ${\mathcal S}(\Omega)$ whose volume is $V$. 

A hypersurface $X\in {\mathcal S}(\Omega; V)$ is a critical point of the area in ${\mathcal S}(\Omega; V)$ if and only if the following two conditions (C1), (C2) are satisfied (cf. \S \ref{s:eigen}). 

(C1) The mean curvature of $X$ is constant.

(C2) $X$ is orthogonal to $\Pi_1\cup\Pi_2$ on the boundary. 

Let $X\in {\mathcal S}(\Omega; V)$. 
Any smooth one-parameter family $X[\epsilon]$ ($-\epsilon_0\!<\!\epsilon\!<\!\epsilon_0, \ \epsilon_0>0$) of hypersurfaces in ${\mathcal S}(\Omega; V)$ satisfying $X[0]=X$ is called a volume-preserving admissible variation of $X$. 

Assume now that a hypersurface $X\in {\mathcal S}(\Omega; V)$ satisfies both (C1) and (C2). 
$X$ is said to be stable if the second variation of the area for all volume-preserving admissible variations of $X$ is nonnegative. 
If $X$ is not stable, then $X$ is said to be unstable.


\begin{figure}[hb]
\begin{center}
\tabcolsep=1mm
		\begin{tabular}{ c c }
			\includegraphics[width=30mm]{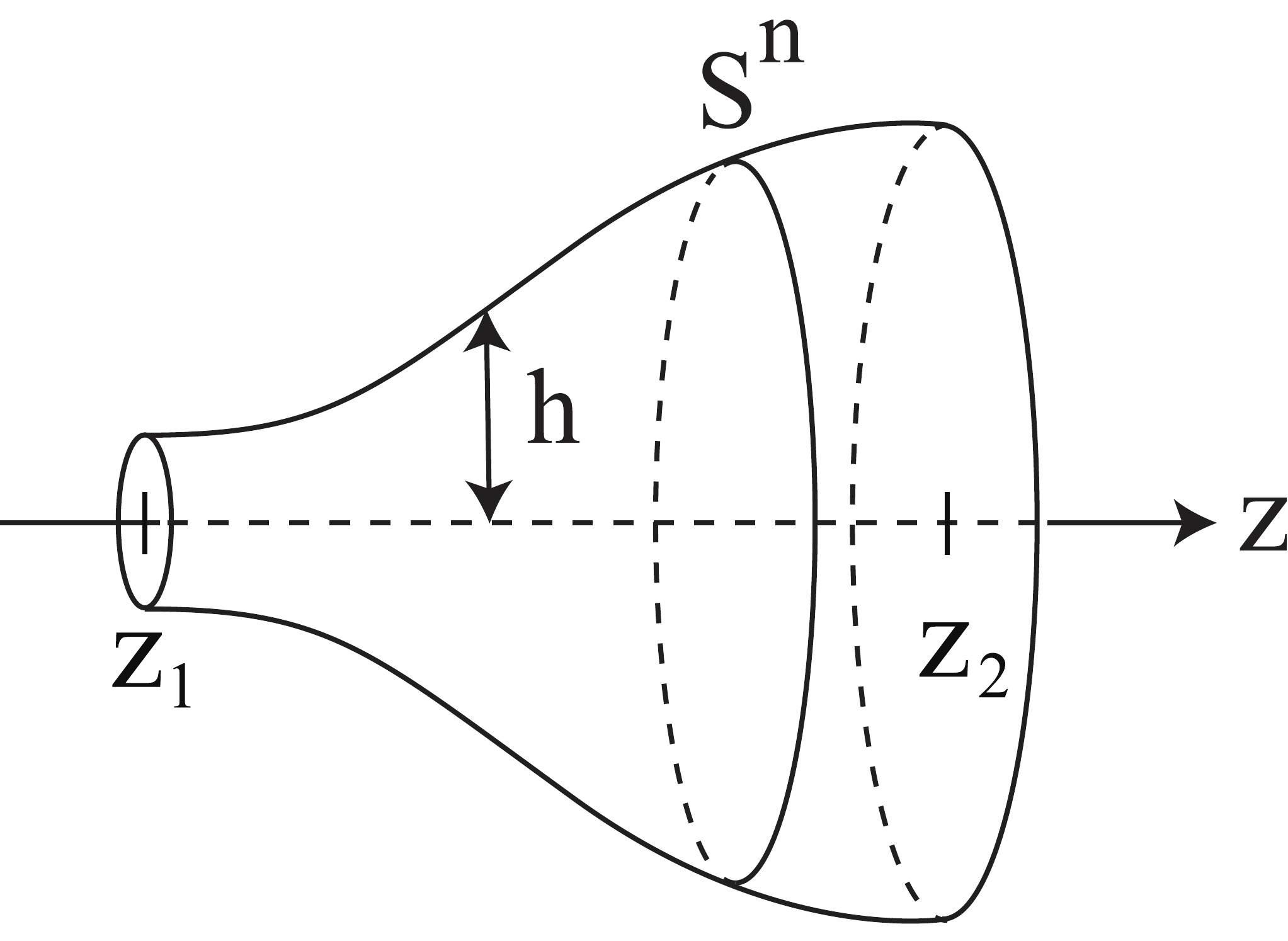} &
			\includegraphics[width=40mm]{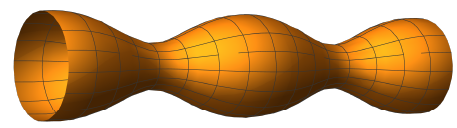} \\
			(a) & (b) \\
		\end{tabular}
\caption{
(a) A hypersurface of revolution in ${\mathbb R}^{n+2}$ with the $z$-axix as the axis of rotation. Its profile curve is given by a function $h=h(z)$ that represents the distance from the $z$-axis. (b) An unduloid.
}
\label{fg:unduloid}
\end{center}
\end{figure}


If $X$ does not have self-intersection (that is, $X$ is an embedded hypersurface), then, by using the Alexandrov reflection method (cf. \cite{K1986}), $X$ is a hypersurface of revolution with the $z$-axis as the rotation axis. 
Denote by $h=h(z)$  ($z_1\le z\le z_2$) (Fig.~\ref{fg:unduloid}(a)) the distance from the $z$-axis to the profile curve of $X$. We sometimes denote $X$ by $X[h]$. 
The mean curvature $H$ of $X$ with respect to the outward-pointing unit normal $\nu$ is
\begin{equation*}
	H
	=
	\frac{1}{n+1} \Big[ \frac{ h_{zz} }{ (1+h_z^2)^{3/2} } - \frac{n}{h\sqrt{1+h_z^2}} \Big].
\end{equation*}
The geometry of the profile curves of the solutions of the ordinary differential equation ``$H=$constant'' is well-known \cite{D, HY}. 
Among them, embedded hypersurfaces satisfying the above (C2) are classified into the following four classes. 

(S1) $(n+1)$-dimensional spheres $S^{n+1}$ included in $\Omega^\circ$,

(S2) hemispheres of $S^{n+1}$ with boundary in either $\Pi_1$ or $\Pi_2$, 

(S3) parts of cylinders with boundary in $\Pi_1 \cup \Pi_2$, and

(S4) $k/2$ period of an unduloid perpendicular to $\Pi_1 \cup \Pi_2$ on the boundary, where $k$ is any positive integer. 

Here the unduloids are two parameter family of hypersurfaces of revolution, which are periodic in the $z$-direction (Fig.~\ref{fg:unduloid}(b)). The shape of its half-period is like Fig.~\ref{fg:unduloid}(a). 

As for the stability, spheres and hemispheres are stable, 
because spheres are the minimizers of the area among all closed hypersurfaces enclosing the same volume, and hemispheres are the minimizers of the area among all compact hypersurfaces enclosing the same volume with non-empty boundary either in $\Pi_1$ or in $\Pi_2$. 
A cylinder with radius $r$ and length $L$ is stable if and only if $\pi r\ge \sqrt{n} L$ (cf. \S \ref{sec:sign}). 
Since any one period of an unduloid from a neck to the next neck is a maximum stable subset of the unduloid for fixed boundary problem \cite{KP2017}, any unduloid with $k/2$ period ($k\ge 2$) is unstable for our free boundary problem. 
Therefore, only the stability of any half period of an unduloid must be studied. 

For each $s \in(0, 1)$, there exists a unique half period of an unduloid (which we denote by ${\mathcal U}(s)$, or simply ${\mathcal U}$) with neck in $\Pi_1$ and bulge in $\Pi_2$ such that 
\begin{equation}\label{nonuni}
	s = 1-\frac{\text{the radius of the neck}}{\text{the radius of the bulge}}
\end{equation}
holds. 
Each limit $s=0, 1$ corresponds to the cylinder, the half sphere, respectively. 

On the stability of ${\mathcal U}$, the following results are known 
\cite{Athanassennas,Vogel,PR,Li}.
\begin{itemize}
\item[{(a)}]
For any $n \geq 1$, if ${\mathcal U}$ is sufficiently close to a hemisphere, then ${\mathcal U}$ is unstable. 
\item[{(b)}]
For $1 \!\leq \! n \!\leq 6$, ${\mathcal U}$ is unstable.
\item[{(c)}]
For $7 \!\leq n\! \leq \! 9$ (resp. $n \!\geq \!10$), if ${\mathcal U}$ is sufficiently close to a cylinder, then ${\mathcal U}$ is unstable (resp. stable).
\item[{(d)}]
For $n \geq 8$, there exists some ${\mathcal U}$ that is stable. 
\end{itemize}

In this paper, we comprehensively examine the stability of unduloids ${\mathcal U}(s)$ in all dimensions as follows. 
Denote by $V(s)$ the volume enclosed by ${\mathcal U}(s)$. 
\begin{itemize}
\item[{(I)}]
For any $n \geq 1$, 
${\mathcal U}(s)$ is stable (resp. unstable)if and only if $V'(s) \leq 0$ (resp.\ $V'(s) > 0$).
\item[{(II)}]
When $1 \leq n \leq 6$, ${\mathcal U}(s)$ is unstable for all $s \in (0,1)$. 
\item[{(III)}]
When $7 \leq n \leq 9$, there exists numbers $s_1$, $s_2$ satisfying $V'(s_1)=V'(s_2)=0$ and $0<s_1<s_2<1$, such that if $s \in [s_1,s_2]$ (resp.\ $s \in (0,s_1) \cup (s_2,1)$), then ${\cal U}(s)$ is stable (resp.\ unstable). 
\item[{(IV)}]
When $n \ge 10$, there exists some $s_2 \in (0,1)$ satisfying $V'(s_2)=0$, such that if $s \in (0,s_2]$ (resp.\ $s \in (s_2,1)$), then ${\cal U}(s)$ is stable (resp.\ unstable).
\end{itemize}
Especially, the existence of a stable unduloid for $n=7$ is found for the first time in this paper.

The stability of the unduloid depending on $s$ is summarized in Table~\ref{tbl:lambda2}, and values of $s_k \; (k=0,1,2,3)$ in Table~\ref{tbl:lambda2} numerically obtained are presented in Table~\ref{tbl:sk}. 
The precise analysis is given in a full paper \cite{KM}. 
Particularly, our original idea is to apply the bifurcation theory constructed in \cite{KPP2017} to judging the stability.

\section{The eigenvalue problem associated with the second variation of the area}\label{s:eigen}

In order to judge the stability of an unduloid, it is sufficient to study only axially-symmetric variations (cf. \cite{KP2006}). 
For this reason, we consider any variation $h[\epsilon](z)=h_0(z)+h_1(z)\epsilon+{\mathcal O}(\epsilon^2)$ of the profile curve $h=h_0(z)$, where $\epsilon$ is the variation parameter. 

Denote by $a_n$ the $n$-dimensional volume of the $n$-dimensional unit sphere. 
The first variations of the area and the volume are
\begin{align}
	\frac{dA(X[h[\epsilon]])}{d\epsilon}\Big|_{\epsilon=0}
	=&
	-(n+1) a_n \int_{z_1}^{z_2} H_0 h_0^n h_1 dz \nonumber\\
	&+
	a_n \Big[
		\frac{h_0^n h_{0z}  }{\sqrt{ 1+h_{0z}^2 }}h_1
	\Big]_{z_1}^{z_2},
\label{A1}
\\
	\frac{dV(X[h[\epsilon]])}{d\epsilon}\Big|_{\epsilon=0}
	=&
	a_n \int_{z_1}^{z_2} h_0^n h_1 dz.
\label{V1}
\end{align}
From \eqref{A1} and \eqref{V1}, one sees that the hypersurface $X[h_0]$ is a critical point of the area for volume-preserving admissible variations  if and only if 
\begin{gather}
	H(z)={\text{const.}},
\;\;\;
 h_{0z}(z_1)=h_{0z}(z_2)=0
\label{eq:h0}
\end{gather}
hold.
Assume now that $h_0$ satisfies \eqref{eq:h0}. 
We assume that $h[\epsilon]$ is a volume-preserving variation. 
Then, the second variation of the area is given by 
\begin{align*}
	\frac{d^2\!A(X[h[\epsilon]])}{d\epsilon^2}\Big|_{\epsilon=0}
	\!\!=\!
	-a_n\int_{z_1}^{z_2}
	h_1 \mathcal{L} h_1 dz 
	\!+\!
	a_n \left[  h_1 \sigma(z) h_{1z} \right]_{z_1}^{z_2},
\end{align*}
where
\begin{align*}
	\mathcal{L}
	:=
	\frac{d}{dz}\left(  \sigma(z) \frac{d}{dz} \right) + \frac{ nh_0^{n-2} }{ \sqrt{1+h_{0z}^2} }, \;\;
	\sigma(z)
	:=
	\frac{ h_0^n }{ \left( 1+h_{0z}^2 \right)^{3/2} }.
\end{align*}
The following eigenvalue problem will be used to analyze the stability of the unduloid. 
\begin{equation} \label{A2}
	{\cal L} \varphi =\! -\lambda \varphi \  (z\! \in\! [z_1,  z_2]), \  
	\varphi_{z}(z_1) = \varphi_{z}(z_2) = 0.
\end{equation}

\section{Outline of the proof 
}

As in \S \ref{intro}, we denote by ${\mathcal U}(s)$ the half period of an unduloid with profile curve $h=h_0(z)$ satisfying \eqref{nonuni}. 
Denote its area, volume, mean curvature, and eigenvalues of the problem \eqref{A2} by $A(s)$, $V(s)$, $H(s)$, and $\lambda_i(s)$, respectively. 
Then, they are continuous functions in $0 \le s <1$, and 
$\lim_{s \to 1-0} f(s) = f(1)$ holds for $f=A, V$, and $H$.

First assume $\lambda_2(s)<0$. Set
$$
u=-\bigl(
\int_{z_1}^{z_2} h_0^n\varphi_2\;dz
\bigr)
\bigl(
\int_{z_1}^{z_2} h_0^n\varphi_1\;dz
\bigr)^{-1}
\varphi_1+\varphi_2,
$$
where $\varphi_j$ is an eigenfunction belonging to $\lambda_j$ ($j=1, 2$). Then, there exists a volume-preserving admissible variation of ${\mathcal U}(s)$ for which the profile curve is $h[\epsilon](z)=h_0(z)+u(z)\epsilon+{\mathcal O}(\epsilon^2)$. The second variation $A_2$ of the area for this variation is 
$$
A_2=a_n\int_{z_1}^{z_2} (\lambda_1\varphi_1^2+\lambda_2\varphi_2^2) \;dz <0.
$$
Hence ${\mathcal U}(s)$ is unstable, that is
\begin{equation}
	\lambda_2 < 0 \;\;\; \Longrightarrow \;\;\; \mbox{unstable.}
\label{cr0}
\end{equation}
For any ${\mathcal U}(s)$, we can show the following inequalities.
\begin{equation}\nonumber
	\lambda_1(s) < 0 < \lambda_3(s), \;\;\; \forall s \in (0,1).
\label{lambda13}
\end{equation}

\subsection{Sign of $\lambda_2$}
\label{sec:sign}

\begin{figure}[hb]
\begin{center}
\tabcolsep=2mm
		\begin{tabular}{ c c }
			\includegraphics[width=35mm]{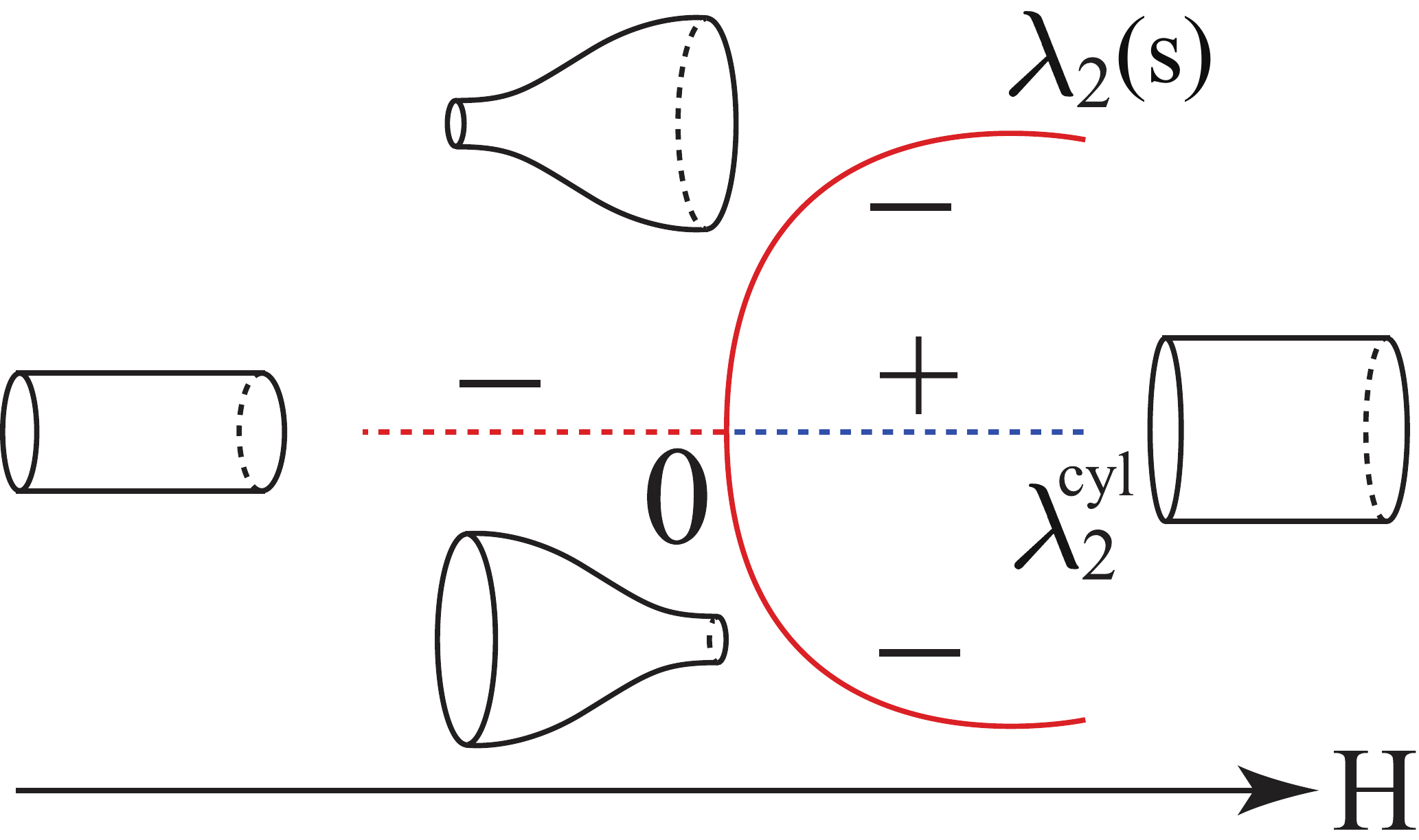} &
			\includegraphics[width=35mm]{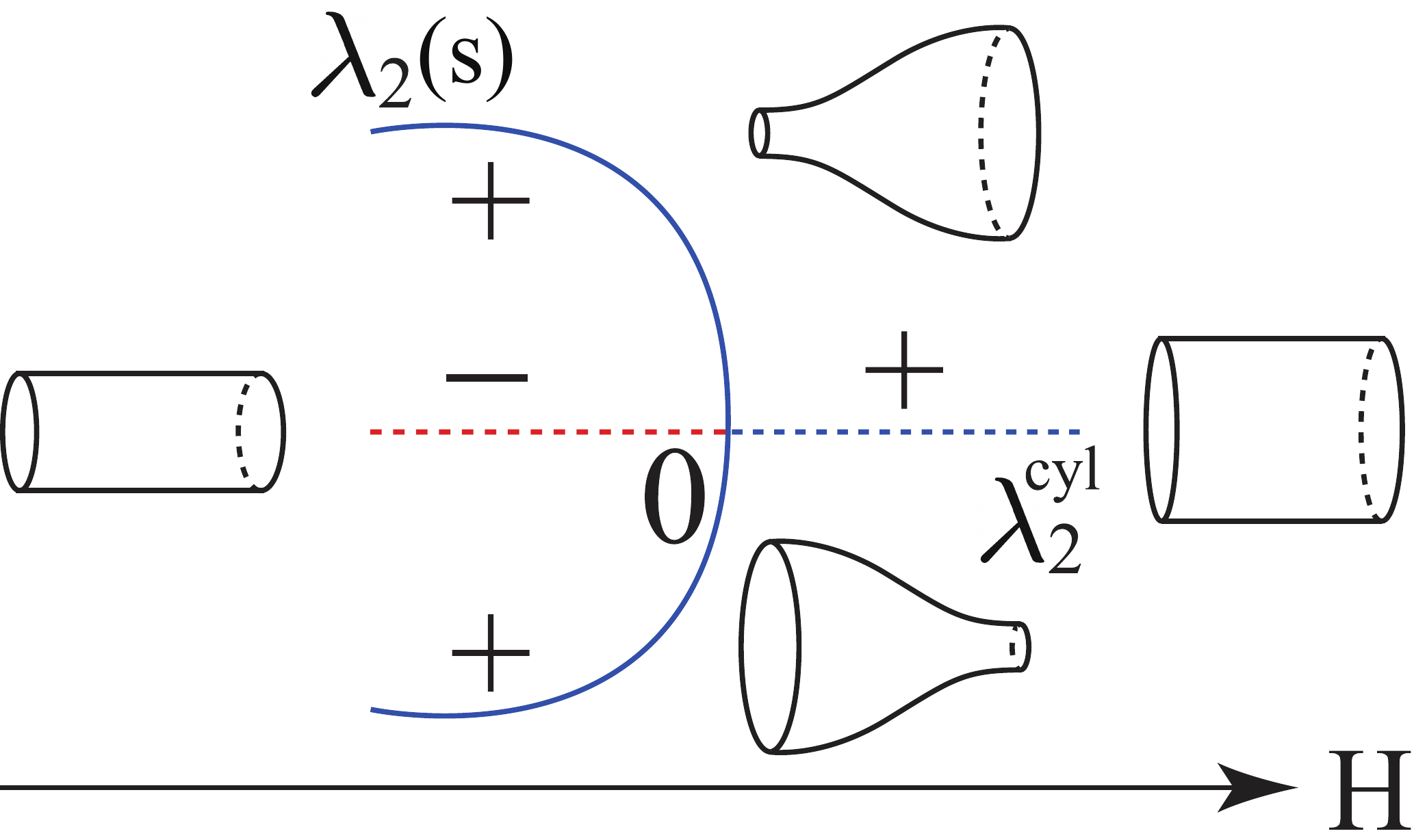} \\
			(a) & (b) \\
		\end{tabular}
\caption{Diagrams representing criterion \eqref{cr1}. If the mean curvature of unduloid is (a) larger (resp.\ (b) smaller) than that of the critical cylinder that is $H(s) > H(0)$ (resp.~$H(s) < H(0)$), $\lambda_2 $ is negative (resp.\ positive). Note that $H$ takes negative values in the present convention.}
\label{fig:bif}
\end{center}
\end{figure}

If one puts $z_1 = 0, z_2  = L \; (>0) $, then the eigenvalue $\lambda_i^{\rm cyl}$ of the cylinder with radius $r$ in ${\mathcal S}(\Omega)$ is 
\begin{equation}
	\lambda_i^{\rm cyl}
	=
	\left(
		\Big[ \frac{ (i-1) \pi r}{L} \Big]^2-n
	\right) r^{n-2},
\;\;\;
	i=1,2,3,\cdots.
\label{lambda^cyl}
\end{equation}
From Eq.~\eqref{lambda^cyl}, one can see that if 
\begin{align} \nonumber
	r < r_c := \frac{\sqrt{n} L}{\pi},
\end{align}
then $\lambda_2^{\rm cyl} < 0$ holds and so, the cylinder with radius $r$ and length $L$ is stable if and only if $r\ge r_c$ holds (cf. \cite{Li}). 
We call the cylinder with critical radius $r_c$ (therefore $\lambda_2^{\rm cyl}=0$) {\it a critical cylinder}.
From the critical cylinder, two branches of unduloid
emanate (see Fig.~\ref{fig:bif}). 
If the mean curvature $H(s)$ of emanating unduloid is larger (resp.\ smaller) than $H(0)$, the second eigenvalue $\lambda_2(s)$ of unduloid is negative (resp.\ positive) \cite{KPP2017}, that is
\begin{align}
	H(s) \gtrless H(0) \;\; \Longrightarrow \;\; \lambda_2 (s) \lessgtr 0
\label{cr1}
\end{align}
for $s$ close to $0$.
This criterion is visualized in Fig.~\ref{fig:bif}.

In order to know when $\lambda_2(s)$ changes its sign, the following criteria are quite useful.
\begin{align}
	\mbox{$H' \neq 0$ at $s$}
	\Longrightarrow 
	\mbox{$\lambda_2$ does not change its sign at $s$},
\label{cr21}
\end{align}
\begin{align}
	&\mbox{$H'=0$ \& $H''\ne 0$ \& $V'\ne 0$ at $s$}\;\;\;
	 \;\;\; \nonumber \\
	&\Longrightarrow \;\;\;
	\;\;\;
	\mbox{$\lambda_2=0$ \& $\lambda_2$ changes its sign at $s$}.
\label{cr22}
\end{align}
The criteria \eqref{cr21} and \eqref{cr22} mean that, under the assumption that $H''(s)\ne 0$ and $V'(s)\ne 0$, $\lambda_2(s)$ changes its sign when $H'(s)$ does. As for the possibility that $\lambda_2(s)$ vanishes, it is proved that 
\begin{align}
	H'(s) \neq 0 \;\;\; \Longrightarrow \;\;\; \lambda_2(s) \neq 0.
\label{cr23} 
\end{align}
Thus, once the sign of $\lambda_2(s)$ near $s=0$ is determined by \eqref{cr1}, the sign of $\lambda_2$ in the full range of $s$ is known by investigating $H''(s)$ and $V'(s)$ at zeros of $H'(s)$.

\subsection{Criteria when $\lambda_2 \geq 0$}
\label{sec:positive}

When $\lambda_2 (s) \geq  0$, by using a modification of Corollary 1.1 in \cite{Ko2}, 
we have the following observations.

{\it When $\lambda_2(s) = H'(s)= 0$ and $V'(s) \ne 0$ hold, ${\cal U}(s)$ is unstable}. Namely, the following holds,
\begin{equation}
	\mbox{ $\lambda_2(s)=H'(s)= 0$ \& $V'(s) \neq 0$}
 \Longrightarrow 
	\mbox{unstable. }
\label{cr3}
\end{equation}
When $\lambda_2(s) > 0$, we have the following criterion of the stability for ${\mathcal U}(s)$.
\begin{align}
\begin{split}
	\lambda_2(s) > 0 \;
	\&
	\begin{cases}
	H'(s) V'(s) < 0 \\
	H'(s) V'(s) \geq 0
	\end{cases}
\Rightarrow 
	\begin{cases}
	\mbox{unstable}\\
	\mbox{stable}
	\end{cases}
\end{split}
\label{cr4}
\end{align}

\begin{table}[htb]
\begin{center}
\begin{minipage}[c]{0.45\textwidth}
\begin{center}
\tabcolsep=1.6 mm

\caption{The sign of $H'(s), \lambda_2(s)$, and $V'(s)$ and the stability of unduloid ${\mathcal U}(s)$ in ${\mathbb R}^{n+2}$ as functions of non-uniformness parameter $s \in (0,1)$ in Class A, B, C, and D.}
\label{tbl:lambda2}

$\mbox{Class A ($ 1\leq n \leq 6$)}$\\
\begin{tabularx}{35mm}{c|  c     }
\Hline
$s$ &   \\
\hline
$H'(s)$ &  $+$   \\
\hline
$\lambda_2(s)$ &  $-$   \\
\hline
$V'(s)$ & +   \\
\hline
stability &  unstable     \\
\Hline
\end{tabularx}
 
\vspace{5pt}

$\mbox{Class B ($n=7$)}$\\
\begin{tabularx}{68mm}{c|  c | c | c | c | c | c | c | c | c }
\Hline
$s$ & & $s_0$ & & $s_1$ & & $s_2$ & & $s_3$ &  \\
\hline
$H'(s)$ &  $+$  & 0 & \multicolumn{5}{c|}{$-$} & $0$ & $+$\\
\hline
$\lambda_2(s)$ &  $-$  & 0 & \multicolumn{5}{c|}{$+$} & $0$ & $-$\\
\hline
$V'(s)$ & \multicolumn{3}{c|}{$+$}  & 0 & $-$& 0 &  \multicolumn{3}{c}{$+$} \\
\hline
stability &  \multicolumn{3}{c|}{unstable}  &  \multicolumn{3}{c|}{stable} & \multicolumn{3}{c}{unstable} \\
\Hline
\end{tabularx}

\vspace{5pt}

$\mbox{Class C ($ 8 \leq n \leq 9$)}$\\
\begin{tabularx}{68mm}{c|  c | c | c | c | c | c | c  }
\Hline
$s$ & & $s_1$ & & $s_2$ & & $s_3$ &    \\
\hline
$H'(s)$ &  \multicolumn{5}{c|}{$-$}  & 0 & $+$ \\
\hline
$\lambda_2(s)$ &  \multicolumn{5}{c|}{$+$}  & 0 & $-$\\
\hline
$V'(s)$ & $+$  & 0 & $-$& 0 &  \multicolumn{3}{c}{$+$} \\
\hline
stability &  unstable  &  \multicolumn{3}{c|}{stable} & \multicolumn{3}{c}{unstable} \\
\Hline
\end{tabularx}

\vspace{5pt}

$\mbox{Class D ($ n \geq 10 $)}$\\
\begin{tabularx}{50mm}{c|  c | c | c | c | c }
\Hline
$s$ & & $s_2$ & & $s_3$ & \\
\hline
$H'(s)$ &  \multicolumn{3}{c|}{$-$}  & 0 & $+$ \\
\hline
$\lambda_2(s)$ &  \multicolumn{3}{c|}{$+$}  & 0 & $-$\\
\hline
$V'(s)$ & $-$  & 0 & \multicolumn{3}{c}{$+$} \\
\hline
stability &  \multicolumn{2}{c|}{stable} & \multicolumn{3}{c}{unstable} \\
\Hline
\end{tabularx}

\end{center}
\end{minipage}
\end{center}
\end{table}

\subsection{Stability of unduloids in ${\mathbb R}^{n+2}$ ($ n \in {\mathbb N} $)}
\label{sec:analysis}

\begin{figure}
\begin{center}
\begin{minipage}[c]{0.45\textwidth}
\begin{center}
\tabcolsep=5mm
		\begin{tabular}{ c }
			{$n=6 \in \mbox{Class A}$}\\
			\includegraphics[scale=0.23]{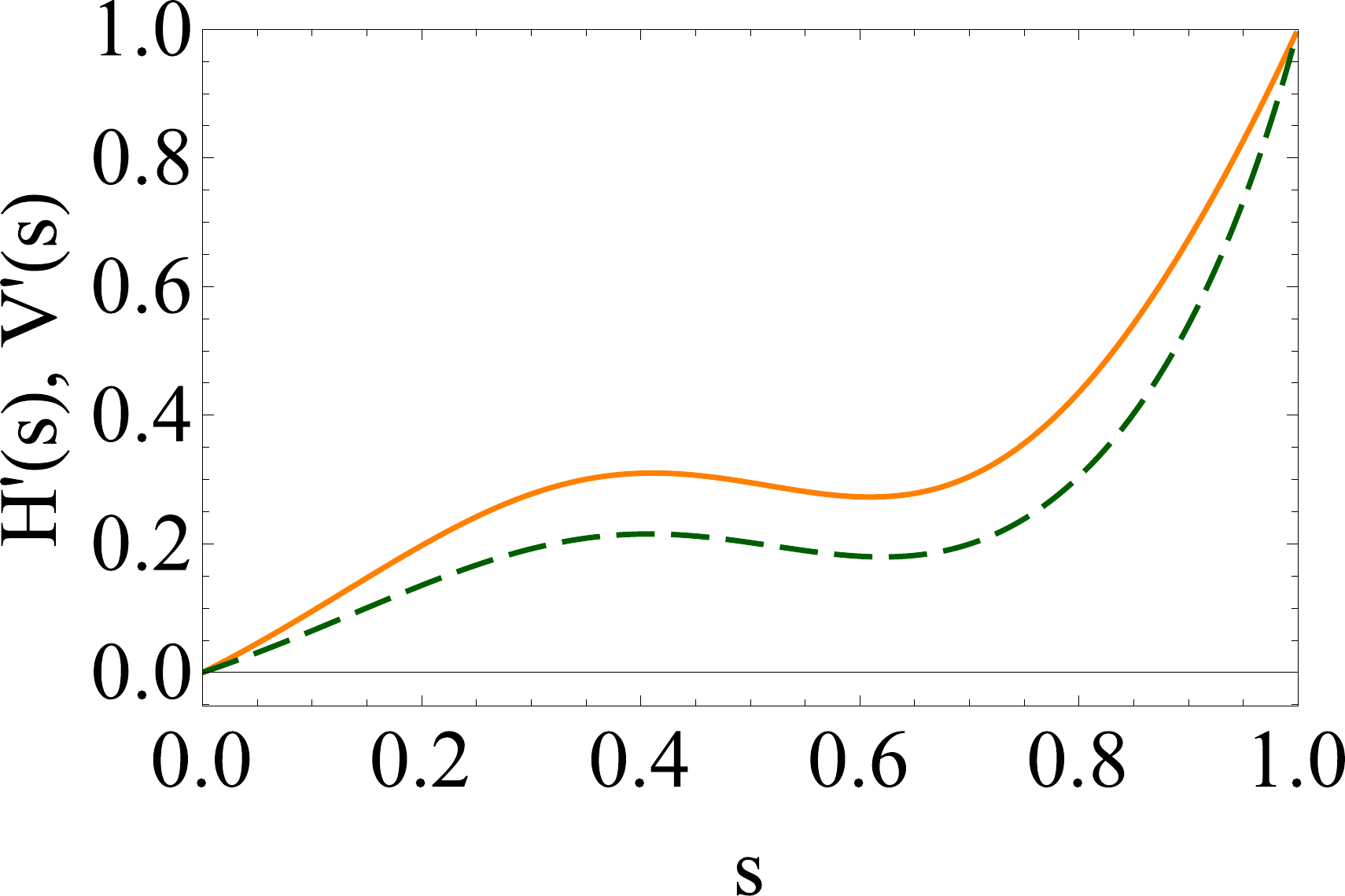} \\
			{$n=7 \in \mbox{Class B}$}\\
			\includegraphics[scale=0.23]{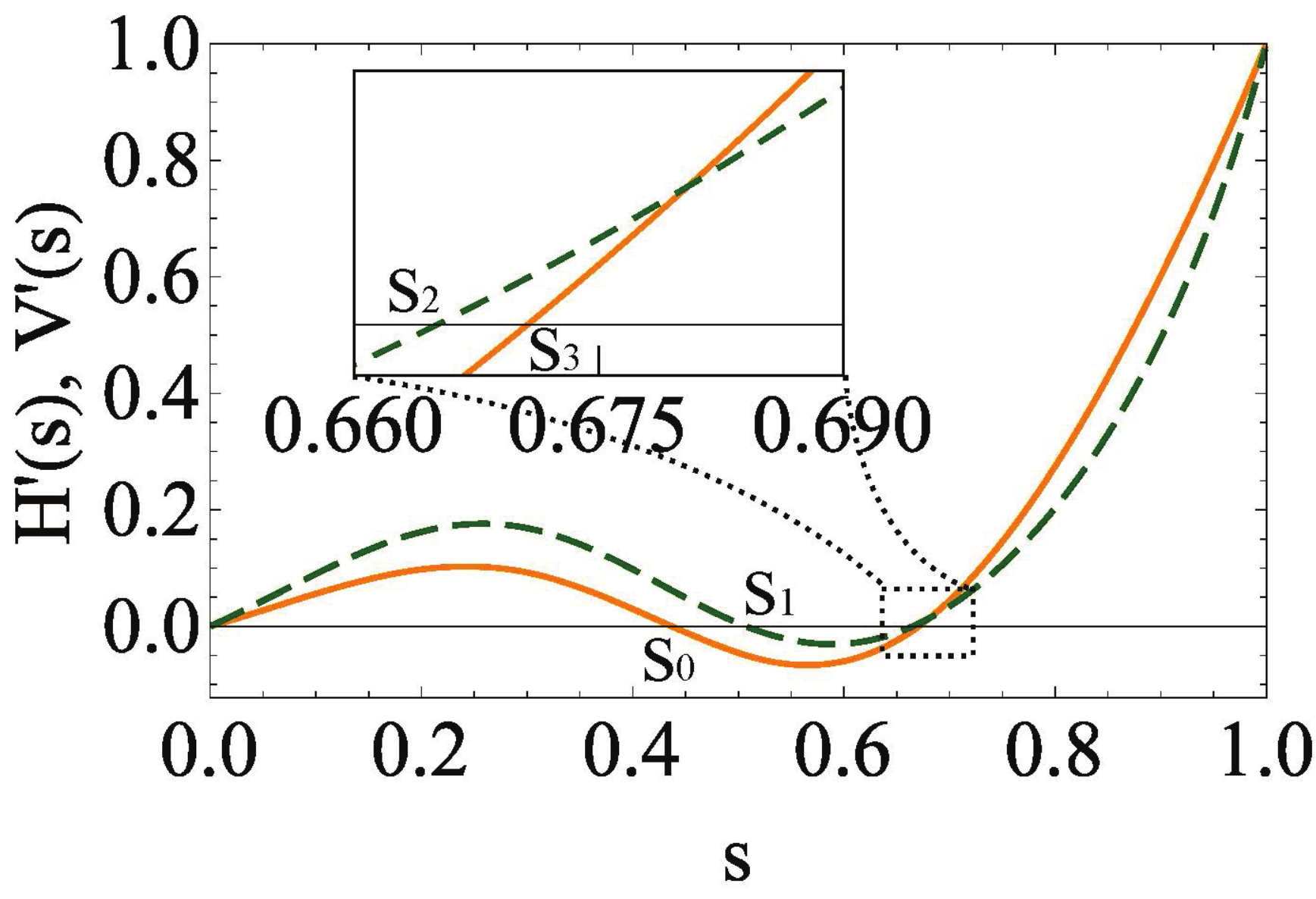} \\
			{$n=8 \in \mbox{Class C}$}\\
			\includegraphics[scale=0.23]{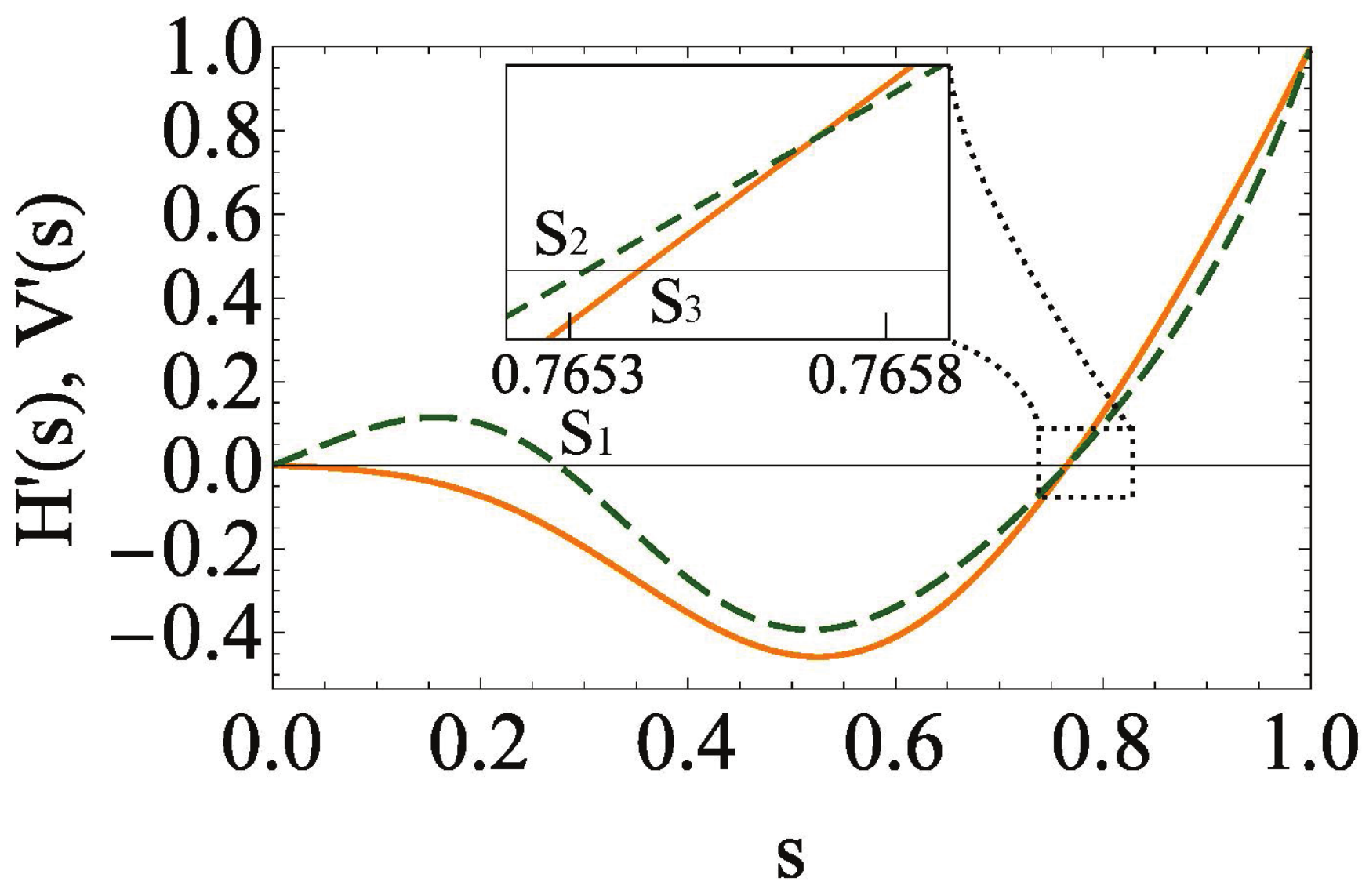} \\
			{$n=10 \in \mbox{Class D}$}\\
			\includegraphics[scale=0.23]{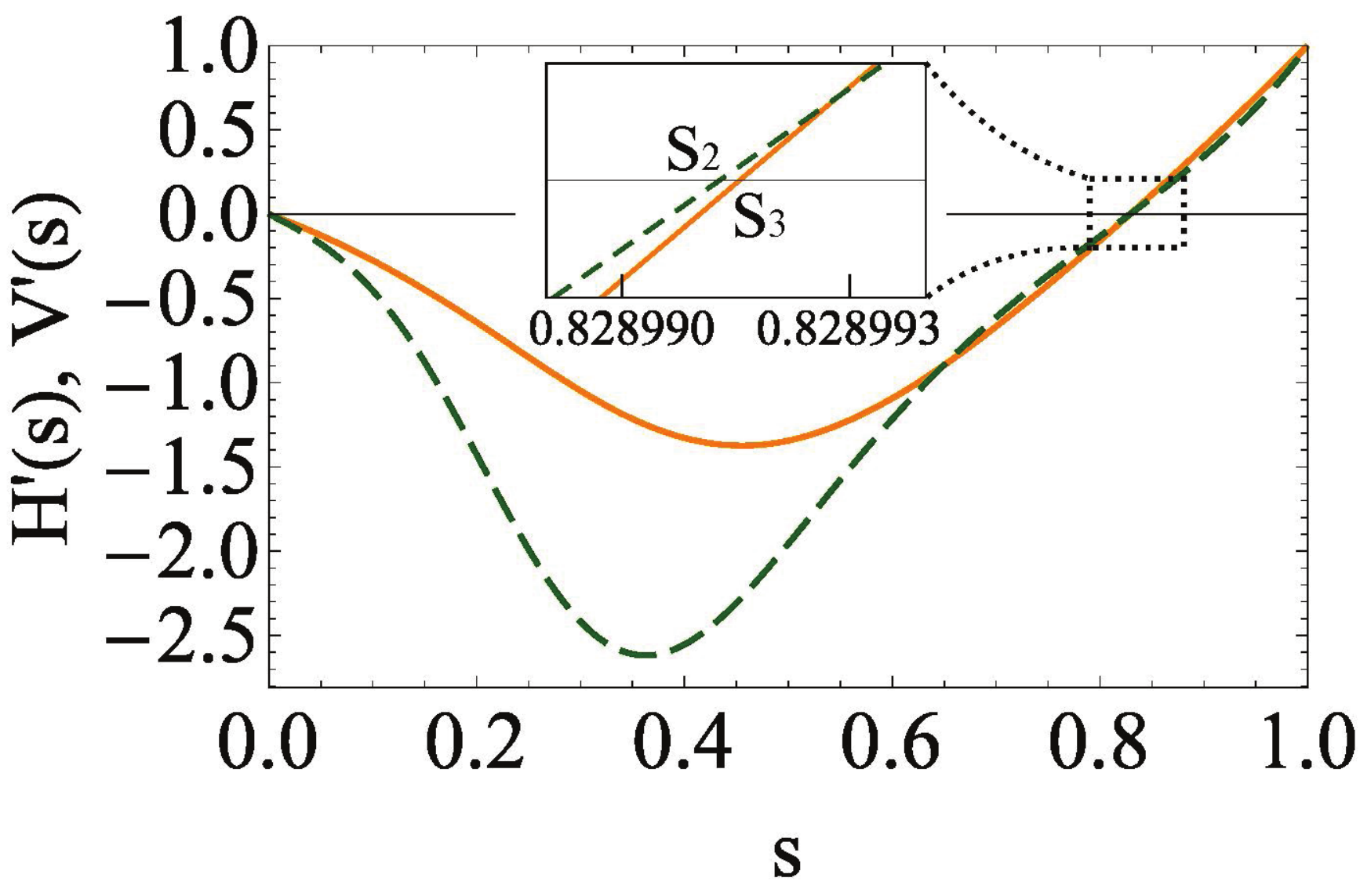} \\
		\end{tabular}
\caption{$H'(s)$ (solid orange line) and $V'(s)$ (dashed green line) of unduloid ${\mathcal U}(s)$ for $n=6, 7, 8$, and $10$ from the top to the bottom.}
\label{fig:diagrams}
\end{center}
\end{minipage}
\end{center}
\end{figure}

We can classify the dimensions into four classes, A ($ 1\leq n \leq 6$), B ($n=7$), C ($ 8 \leq n \leq 9$), and D ($ n \geq 10 $) as in Table~\ref{tbl:lambda2}.  
The numerical plots of $H'(s)$ and $V'(s)$
are presented in Fig.~\ref{fig:diagrams}.
There $H'(s)$ and $V'(s)$ are normalized by $\lim_{s \to 1-0}| H'(s) | $ and $\lim_{s \to 1-0} | V'(s) |$, respectively.
By using this with 
the stability criteria \eqref{cr0} and \eqref{cr1}--\eqref{cr4}, one can show the stability results (I)-(IV) for all unduloids stated in \S \ref{intro}. 
In this paper, we give a detailed proof only for Class B ($n=7$) as follows. 
In Fig.~\ref{fig:diagrams}, we observe that both $H'(s)$ and $V'(s)$ have two simple zeros, which we denote by $s_k$ ($ k=0,1,2,3$) as
\begin{align}\nonumber
\begin{split}
	H'(s_0)= V'(s_1) & =V'(s_2)=H'(s_3)=0,  \\
	0<s_0<s_1  &  <s_2<s_3<1.
\end{split}
\label{n7s}
\end{align}
From the behavior of $H'(s)$, one knows that $\lambda_2(s)$ vanishes and changes its sign only at $s=s_0$ and $s=s_3$ with criteria \eqref{cr21}-\eqref{cr23}. From this and the behavior of $H'(s)$ with criterion \eqref{cr1}, we see that $\lambda_2(s) < 0 $ (resp.\ $\lambda_2(s) \geq 0 $) for $s \in (0,s_0) \cup (s_3,1)$ (resp.\ $s \in [s_0,s_3]$). Hence, ${\mathcal U}(s)$ for $s \in (0,s_0) \cup (s_3,1)$ is unstable with criterion \eqref{cr0}. Since $\lambda_2(s) \geq 0 $ for $s \in [s_0,s_3]$, we have to see also the behavior of $V'(s)$ in order to use criteria~\eqref{cr3} and \eqref{cr4}. From Fig.~\ref{fig:diagrams}, $V'(s)$ vanishes at neither $s=s_0$ nor $s=s_3$, with which \eqref{cr3} implies that ${\mathcal U}(s_0)$ and ${\mathcal U}(s_3)$ are unstable. Since $H'(s)V'(s) < 0$ (resp.\ $ H'(s)V'(s) \geq 0 $), ${\mathcal U}(s)$ is unstable (resp.\ stable) for $ s \in (s_0,s_1) \cup (s_2,s_3) $ (resp.\ $ s \in [ s_1, s_2 ]$).

\begin{table}
\begin{center}
\begin{minipage}[c]{0.45\textwidth}
\begin{center}
\tabcolsep=4 mm

\caption{Values of $s_k\; (k=0,1,2,3)$ for several $n$.
}
\label{tbl:sk}

\vspace{5pt}

\begin{tabularx}{75mm}{c  c  c  c }
\Hline
$n$ & 7 & 8 & 9 \\
\hline
$s_0$ &  0.437 & n/a & n/a \\
\hline
$s_1$ &  0.507 & 0.275 03 & 0.093 270 8 \\
\hline
$s_2$ &  0.665 & 0.765 33 & 0.803 961 7  \\
\hline
$s_3$ &  0.671 & 0.765 41 & 0.803 966 2  \\
\Hline
\end{tabularx}

\vspace{10pt}

\begin{tabularx}{70mm}{c  c  c }
\Hline
$n$ &  10   & 11 \\
\hline
$s_0$ &   n/a & n/a  \\
\hline
$s_1$ &   n/a  & n/a \\
\hline
$s_2$ &   0.828 991 30 & 0.847 468 517 \\
\hline
$s_3$ &   0.828 991 56  & 0.847 468 533 \\
\Hline
\end{tabularx}

\end{center}
\end{minipage}
\end{center}
\end{table}

%

\acknowledgments

The authors are partially supported by JST CREST Grant Number JPMJCR1911. 
This work was also supported in part by JSPS KAKENHI Grant Numbers 
JP18H04487, JP20H04642, JP20H01801 (K.M.),
JP18K03652,
and
JP22K03623 (U.M.).

\references


\begin{thebibliography}{11}


%
\bibitem{K1986}
M.\ Koiso, Symmetry of hypersurfaces of constant mean curvature with symmetric boundary, 
Math. Zeit., {\bf 191} (1986), 567--574.
%
\bibitem{D}
C.\ Delaunay, 
Sur la surface de r\'evolution dont la courbure mayenne est constante, 
J. Math.
Pures. Appl., S\'er. 1, {\bf 6} (1841), 309--320.
%
\bibitem{HY}
W.-Y.\ Hsiang and W.-C.\ Yu,
A generalization of a theorem of Delaunay, 
J. Diff. Geom., {\bf 16} (1981), 161--177.
%
\bibitem{KP2017}
M.\ Koiso and B.\ Palmer, 
Higher order variations of constant mean curvature surfaces, 
Calc. Var. Partial Differential Equations, {\bf 56} (2017), 159. 
%
\bibitem{Li}
H.\ Li, Y.\ Xia and C.\ Xiong,
Stability of unduloid bridges with free boundary in a Euclidean slab, 
Science China Mathematics, {\bf 61} (2018), 917--928.
%
\bibitem{Athanassennas}
M.\ Athanassennas,
A variational problem for constant mean curvature surfaces with free boundary, 
J.\ Reine Angew.\ Math.,\ {\bf 377} (1987), 97--107.
%
\bibitem{Vogel}
T.~I.~Vogel,
Stability of a liquid drop trapped between two parallel planes, 
SIAM J.\ Appl.\ Math.,\ {\bf 47} (1987), 516--525.
%
\bibitem{PR}
H.\ L.\ Pedrosa and M.\ Ritor\'e, 
Isoperimetric domains in the Riemannian product of a circle with a simply connected space form and applications to free boundary problems, 
Indiana U.\ Math.\ J.,\ {\bf 48} (1999), 1357--1394. 
%
\bibitem{KM}
M.\ Koiso and U.\ Miyamoto, 
Stability of hypersurfaces with constant mean curvature trapped between two parallel hyperplanes, accepted for publication in JJIAM, arXiv:1905.01705 [math-ph].


\bibitem{KPP2017}
M.~Koiso, B.~Palmer and P.~Piccione, 
Stability and bifurcation for surfaces with constant mean curvature, 
Journal of the Mathematical Society of Japan, {\bf 69} (2017), 1519--1554.
%
\bibitem{KP2006} M.\ Koiso and B.\ Palmer, 
Stability of anisotropic capillary surfaces between two parallel planes, 
Calc. Var. Partial Differential Equations,
{\bf 25} (2006), 275--298.
%
\bibitem{Ko2}
M. Koiso, Deformation and stability of surfaces with constant mean curvature, 
Tohoku Math. J.(2), \textbf{54}\! (2002), \!145--159.
%
\end{thebibliography}
\end{document}